\documentclass[11pt,twoside]{article}
\usepackage[T1]{fontenc}
\usepackage{amsmath,amsfonts,amssymb,amsthm,fullpage,bbm}
\usepackage{fancyhdr,graphicx,color}
\usepackage{epsfig}

\newtheorem{thm}{Theorem}

\theoremstyle{definition}

\theoremstyle{remark}
\newtheorem{rem}{Remark}

\DeclareMathOperator{\cyl}{cyl}
\DeclareMathOperator{\card}{card}

\DeclareMathOperator{\hyp}{hyp}

\DeclareMathOperator{\flow}{flow}

\newcommand{\eps}{\varepsilon}

\def\PP{\mathbb{P}}
\def\RR{\mathbb{R}}
\def\EE{\mathbb{E}}
\def\NN{\mathbb{N}}
\def\ZZ{\mathbb{Z}}

\def\H{\mathcal{H}}

\def\V{\mathcal{V}}

\def\F{\mathcal{F}}

\title{\huge Upper large deviations for maximal flows through a tilted cylinder} 
\author{\Large Marie THERET}
\date{}

\begin{document}
\maketitle
\begin{center}
\vskip-0.8cm
D.M.A., E.N.S.\\
45 rue d'Ulm\\
75230 Paris Cedex 05\\
marie.theret@ens.fr\\
\end{center}

We consider the standard first passage percolation model in $\ZZ^d$ for
$d\geq 2$ and we study the maximal flow from the upper half part to the
lower half part (respectively from the top to the bottom) of a cylinder
whose basis is a hyperrectangle of sidelength proportional to $n$ and
whose height is $h(n)$ for a certain height function $h$. We denote
this maximal flow by $\tau_n$ (respectively $\phi_n$). We emphasize
the fact that the cylinder may be tilted. We look at the probability that
these flows, rescaled by the surface of the basis of the cylinder, are
greater than $\nu(\vec{v})+\eps$ for some positive $\eps$, where
$\nu(\vec{v})$ is the almost sure limit of the rescaled variable $\tau_n$
when $n$ goes to infinity. On one hand, we prove that the speed of decay of
this probability in the case of the variable $\tau_n$ depends on the tail
of the distribution of the capacities of the edges: it can decays
exponentially fast with $n^{d-1}$, or with $n^{d-1} \min(n,h(n))$, or at an
intermediate regime. On the
other hand, we prove that this probability in the case of the variable
$\phi_n$ decays exponentially fast with the volume of the cylinder as soon
as the law of the capacity of the edges admits one exponential moment; the
importance of this result is however limited by the fact that
$\nu(\vec{v})$ is not in general the almost sure limit of the rescaled
maximal flow $\phi_n$, but it is the case at least when the height $h(n)$
of the cylinder is negligible compared to $n$.\\


%%%%%%%%%%%%%%%%%%%%%%%%%%%%%%%%%%%%%%%%%%%

\section{Definitions and main results}

Let $d\geq2$. We consider the graph $(\mathbb{Z}^{d}  ,
\mathbb E ^{d}  )$ having for vertices $\mathbb Z ^{d}  $ and for edges
$\mathbb E ^{d}$, the set of pairs of nearest neighbours for the standard
$L^{1}$ norm. With each edge $e$ in $\mathbb{E}^{d}$ we associate a random
variable $t(e)$ with values in $\mathbb{R}^{+}$. We suppose that the family
$(t(e), e \in \mathbb{E}^{d})$ is independent and identically distributed,
with a common distribution function $F$: this is the standard model of
first passage percolation on the graph $(\mathbb{Z}^d,
\mathbb{E}^d)$. We interpret $t(e)$ as the capacity of the edge $e$; it
means that $t(e)$ is the maximal amount of fluid that can go through the
edge $e$ per unit of time.

The maximal flow $\phi(F_1 \rightarrow F_2 \textrm{ in } C)$ from
$F_1$ to $F_2$ in $C$, for $C \subset \mathbb{R}^d$ (or by commodity the
corresponding graph $C\cap \mathbb{Z}^d$) can be defined properly this way. We will say that an edge 
$e=\langle x,y\rangle$ belongs to a subset $A$ of $\mathbb{R}^{d}$, which
we denote by $e\in A$, if the segment joining $x$ to $y$ (eventually
excluding these points) is included in $A$. We define
$\widetilde{\mathbb{E}}^{d}$ as the set of all the oriented edges, i.e.,
an element $\widetilde{e}$ in $\widetilde{\mathbb{E}}^{d}$ is an ordered
pair of vertices which are nearest neighbours. We denote an element $\widetilde{e} \in \widetilde{\mathbb{E}}^{d}$ by $\langle \langle x,y \rangle \rangle$, where $x$, $y \in \mathbb{Z}^{d}$ are the endpoints of $\widetilde{e}$ and the edge is oriented from $x$ towards $y$. We consider the set $\mathcal{S}$ of all pairs of functions $(g,o)$, with $g:\mathbb{E}^{d} \rightarrow \mathbb{R}^{+}$ and $o:\mathbb{E}^{d} \rightarrow \widetilde{\mathbb{E}}^{d}$ such that $o(\langle x,y\rangle ) \in \{ \langle \langle x,y\rangle \rangle , \langle \langle y,x \rangle \rangle \}$, satisfying:
\begin{itemize}
\item for each edge $e$ in $C$ we have
$$0 \,\leq\, g(e) \,\leq\, t(e) \,,$$
\item for each vertex $v$ in $C \smallsetminus (F_1\cup F_2)$ we have
$$ \sum_{e\in C\,:\, o(e)=\langle\langle v,\cdot \rangle \rangle}
  g(e) \,=\, \sum_{e\in C\,:\, o(e)=\langle\langle \cdot ,v \rangle
    \rangle} g(e) \,, $$
\end{itemize}
where the notation $o(e) = \langle\langle v,. \rangle \rangle$
(respectively $o(e) = \langle\langle .,v \rangle \rangle$) means that there
exists $y \in \mathbb{Z}^d$ such that $e = \langle v,y \rangle$ and $o(e) = \langle\langle v,y \rangle \rangle$ (respectively $o(e) = \langle\langle y,v \rangle \rangle$).
A couple $(g,o) \in \mathcal{S}$ is a possible stream in $C$ from
$F_1$ to $F_2$: $g(e)$ is the amount of fluid that goes through the edge $e$, and $o(e)$ gives the direction in which the fluid goes through $e$. The two conditions on $(g,o)$ express only the fact that the amount of fluid that can go through an edge is bounded by its capacity, and that there is no loss of fluid in the graph. With each possible stream we associate the corresponding flow
$$ \flow (g,o) \,=\, \sum_{ u \in F_2 \,,\,  v \notin C \,:\, \langle
  u,v\rangle \in \mathbb{E}_n^{d}} g(\langle u,v\rangle) \mathbbm{1}_{o(\langle u,v\rangle) = \langle\langle u,v \rangle\rangle} - g(\langle u,v\rangle) \mathbbm{1}_{o(\langle u,v\rangle) = \langle\langle v,u \rangle\rangle} \,. $$
This is the amount of fluid that crosses $C$ from $F_1$
  to $F_2$ if the fluid respects the stream $(g,o)$. The maximal flow through
  $C$ from $F_1$ to $F_2$ is the supremum of this quantity over all possible choices of streams
$$ \phi (F_1 \rightarrow F_2 \textrm{ in }C)  \,=\, \sup \{ flow (g,o)\,|\,
  (g,o) \in \mathcal{S} \}  \,.$$

The maximal flow
$\phi (F_1\rightarrow F_2 \textrm{ in } C)$ can be expressed differently
thanks to the max-flow min-cut theorem (see \cite{Bollobas}). We need some
definitions to state this result.
A path on the graph $\mathbb{Z}^{d}$ from $v_{0}$ to $v_{m}$ is a sequence $(v_{0}, e_{1}, v_{1},..., e_{m}, v_{m})$ of vertices $v_{0},..., v_{m}$ alternating with edges $e_{1},..., e_{m}$ such that $v_{i-1}$ and $v_{i}$ are neighbours in the graph, joined by the edge $e_{i}$, for $i$ in $\{1,..., m\}$.
A set $E$ of edges in $C$ is said to cut $F_1$ from $F_2$ in
$C$ if there is no path from $F_1$ to $F_2$ in $C \smallsetminus
E$. We call $E$ an $(F_1,F_2)$-cut if $E$ cuts $F_1$ from $F_2$ in $C$
and if no proper subset of $E$ does. With each set $E$ of edges we
associate its capacity which is the variable
$$ V(E)\, = \, \sum_{e\in E} t(e) \, .$$
The max-flow min-cut theorem states that
$$ \phi(F_1\rightarrow F_2 \textrm{ in } C) \, = \, \min \{ \, V(E) \, | \, E
\textrm{ is a } (F_1,F_2)\textrm{-cut} \, \} \, .$$

We need now some geometric definitions. For a subset $X$ of
$\mathbb{R}^d$, we denote by $\mathcal{H}^s (X)$ the $s$-dimensional
Hausdorff measure of $X$ (we will use $s=d-1$ and $s=d-2$). The
$r$-neighbourhood $\V(X,r)$ of $X$ for the Euclidean distance $d$ is defined by
$$ \V (X,r) \,=\, \{ y\in \RR^d\,|\, d(y,X)<r\}\,.  $$
If $X$ is a subset of $\RR^d$ included in an hyperplane of $\RR^d$ and of
co-dimension $1$ (for example a non degenerate hyperrectangle), we denote by
$\hyp(X)$ the hyperplane spanned by $X$, and we denote by $\cyl(X, h)$ the
cylinder of basis $X$ and of height $2h$ defined by
$$ \cyl (X,h) \,=\, \{x+t \vec{v} \,|\, x\in X \,,\,  t\in
[-h,h]    \}\,,$$
where $\vec{v}$ is one of the two unit vectors orthogonal to $\hyp(X)$.

Let $A$ be a non degenerate hyperrectangle,
i.e., a box of dimension $d-1$ in $\mathbb{R}^d$. All hyperrectangles will be
supposed to be closed in $\mathbb{R}^d$. We denote by
$\vec{v}$ one of 
the two unit vectors orthogonal to $\hyp (A)$. For $h$ a
positive real number, we consider the cylinder $\cyl(A,h)$.
The set $\cyl(A,h) \smallsetminus \hyp (A)$ has two connected
components, which we denote by $\mathcal{C}_1(A,h)$ and
$\mathcal{C}_2(A,h)$. For $i=1,2$, let $A_i^h$ be
the set of the points in $\mathcal{C}_i(A,h) \cap \mathbb{Z}_n^d$ which have
a nearest neighbour in $\mathbb{Z}^d \smallsetminus \cyl(A,h)$:
$$ A_i^h\,=\,\{x\in \mathcal{C}_i(A,h) \cap
\mathbb{Z}^d \,|\, \exists y \in \mathbb{Z}^d \smallsetminus \cyl(A,h)
\,,\, \langle x,y \rangle \in \EE^d\}\,.$$
Let $T(A,h)$ (respectively $B(A,h)$) be the top
(respectively the bottom) of $\cyl(A,h)$, i.e.,
$$ T(A,h) \,=\, \{ x\in \cyl(A,h) \,|\, \exists y\notin \cyl(A,h)\,,\,\,
\langle x,y\rangle \in \mathbb{E}^d \textrm{ and }\langle x,y\rangle
\textrm{ intersects } A+h\vec{v}  \}  $$
and
$$  B(A,h) \,=\, \{ x\in \cyl(A,h) \,|\, \exists y\notin \cyl(A,h)\,,\,\,
\langle x,y\rangle \in \mathbb{E}^d \textrm{ and } \langle x,y\rangle
\textrm{ intersects } A-h\vec{v}  \} \,.$$
For a given realization $(t(e),e\in \mathbb{E}^{d})$ we define the variable
$\tau (A,h) = \tau(\cyl(A,h), \vec{v})$ by
$$ \tau(A,h) \,=\,  \tau(\cyl(A,h), \vec{v})\,=\, \phi (A_1^h \rightarrow A_2^h
\textrm{ in } \cyl(A,h)) \,,$$
and the variable $\phi(A,h)= \phi(\cyl(A,h), \vec{v})$ by
$$ \phi(A,h) \,=\,\phi(\cyl(A,h), \vec{v}) \,=\, \phi (B(A,h) \rightarrow T(A,h)
\textrm{ in }   \cyl(A,h))\,, $$ 
where $\phi(F_1 \rightarrow F_2 \textrm{ in } C)$ is defined previously.

There exist laws of large numbers concerning these two variables. We
summarize the results here. The law of large numbers for $\tau$ and for
$\phi$ in flat cylinders is the following:
\begin{thm}[Rossignol and Th\'eret \cite{RossignolTheret08b}]
We suppose that
$$ \int_{[0,+\infty[} x \,dF(x) \,<\,\infty\,. $$
Then for every unit vector $\vec{v}$, there exists a constant $\nu(\vec{v})
= \nu(\vec{v}, d,F)$ such that for every non degenerate
hyperrectangle $A$ orthogonal to $\vec{v}$, for every function $h:\NN
\rightarrow \RR^+$ satisfying $\lim_{n\rightarrow \infty} h(n) = +\infty$,
we have
$$ \lim_{n\rightarrow \infty} \frac{\tau(nA,h(n))}{\H^{d-1}(nA)} \,=\,
\nu(\vec{v}) \quad \textrm{in }L^1\,. $$
Moreover, if $0\in A$, where $0$ is the origin of the graph, or if 
$$ \int_{[0,+\infty[} x^{1+\frac{1}{d-1}} \,dF(x) \,<\,\infty\,, $$
then
$$ \lim_{n\rightarrow \infty} \frac{\tau(nA,h(n))}{\H^{d-1}(nA)} \,=\,
\nu(\vec{v}) \quad \textrm{a.s.} $$
If $\lim_{n\rightarrow \infty} h(n)/n
=0$, the same convergences (in $L^1$ and a.s.) hold for $\phi(nA,h(n))$
under the same hypotheses.
\end{thm}
Thanks to the works of Kesten \cite{Kesten:flows} and Zhang \cite{Zhang},
we know that $\nu(\vec{v})>0$ if and only if $F(0)<1-p_c(d)$, where
$p_c(d)$ is the critical parameter for the edge percolation on $\ZZ^d$.
Kesten, Zhang, and finally Rossignol and Th\'eret have proved a law of
large numbers for the variable $\phi(A, h)$ in straight cylinders,
i.e., when $A$ is of the form $\prod_{i=1}^{d-1}[0,k_i] \times
\{0\}$ with $k_i>0$ for all $i=1,...,d-1$, for large $A$ and $h$. Kesten and Zhang have worked in the general
case where the dimensions of the cylinder go to infinity with possibly
different speed. We
present here the result stated by Rossignol and Th\'eret in
\cite{RossignolTheret08b}, with the best conditions on the moment of $F$
and on the height function $h$, but in the more restrictive case where the
cylinder we consider is simply $\cyl(nA, h(n))$:
\begin{thm}[Rossignol and Th\'eret \cite{RossignolTheret08b}]
We suppose that
$$ \int_{[0,+\infty[} x \,dF(x) \,<\,\infty\,. $$
Let $\vec{v}_0=(0,...,0,1)$. For every hyperrectangle $A$ of the form
$\prod_{i=1}^{d-1}[0,k_i]\times \{0\}$ with $k_i>0$ for all $i=1,..,d-1$,
and for every function $h:\NN \rightarrow \RR^+$ satisfying $\lim_{n\rightarrow
\infty} h(n) =+\infty$ and $\lim_{n\rightarrow \infty} \log h(n) / n^{d-1}
=0$, we have
$$ \lim_{n\rightarrow \infty} \frac{\phi(nA,h(n))}{\H^{d-1}(nA)} \,=\,
\nu(\vec{v}) \quad \textrm{a.s. and in }L^1\,. $$
\end{thm}

We investigate the upper large deviations of the variables $\phi$ and
$\tau$. We will prove the following theorem concerning $\tau$:
\begin{thm}
\label{chapitre3thmtau}
Let $A$ be a non degenerate hyperrectangle, and $\vec{v}$ one of the two
unit vectors normal to $A$. Let $h: \mathbb{N} \rightarrow \mathbb{R}^+$
be a height function satisfying $\lim_{n \rightarrow  \infty} h(n) = +\infty$. The upper large
deviations of $\tau(nA,h(n))/\H^{d-1}(nA)$ depend on the
tail of the distribution of the capacities. Indeed, we obtain that:
\newline
i) if the law of the capacity of the edges has bounded support, then for every
  $\lambda >\nu(\vec{v})$ we have
\begin{equation}
\label{chapitre3casborne}
 \liminf_{n\rightarrow \infty} \frac{-1}{\mathcal{H}^{d-1}(nA) \min(h(n), n)} \log
\mathbb{P} \left[ \frac{\tau(nA, h(n))}{\mathcal{H}^{d-1}(nA)} \geq \lambda
  \right]  \,>\,0\,;
\end{equation}
the upper large deviations are then of volume order for height functions
$h$ such that $h(n)/n$ is bounded, and of order $n^d$ if
$\lim_{n\rightarrow \infty} h(n)/n = +\infty$.
\newline
ii) if the capacity of the edges follows the exponential law of parameter
  $1$, then there exists $n_0 (d, A,h)$, and for
  every $\lambda >\nu(\vec{v})$ there exists a positive constant $D$
  depending only on $d$ and $\lambda$ such that for all $n\geq n_0$ we have
\begin{equation}
\label{chapitre3casexp}
\mathbb{P} \left[ \tau(nA, h(n)) \geq \lambda \H^{d-1}(nA)
  \right] \,\geq\, \exp (- D \H^{d-1}(nA))\,.
\end{equation}
\newline
iii) if the law of the capacity of the
  edges satisfies
$$ \forall \theta >0 \qquad \int_{[0,+\infty[} e^{\theta x} dF(x) \,<\,\infty \,,$$
then for all $\lambda > \nu(\vec{v})$ we have
\begin{equation}
\label{chapitre3castsmts}
 \lim_{n\rightarrow \infty} \frac{1}{\mathcal{H}^{d-1}(nA)} \log
\mathbb{P} \left[ \frac{\tau(nA, h(n))}{\mathcal{H}^{d-1}(nA)} \geq \lambda
  \right] \,=\, -\infty \,.
\end{equation}
\end{thm}
We also prove the following partial result concerning the variable $\phi$:
\begin{thm}
\label{chapitre3thmphi}
Let $A$ be a non degenerate
hyperrectangle in $\RR^d$, of normal unit vector $\vec{v}$, and $h: \NN
\rightarrow \RR^+$ be a function satisfying $\lim_{n\rightarrow
  \infty} h(n) = +\infty$. We suppose that the law of the capacities of the
edges admits an exponential moment:
$$ \exists \gamma >0 \qquad \int_{[0,+\infty[} e^{\gamma x} \,dF(x) \,<\,
\infty\,.  $$
Then for every $\lambda > \nu(\vec{v})$, we have 
$$\liminf_{n\rightarrow \infty} \frac{-1}{\H^{d-1}(nA)h(n)}\log \PP [\phi (nA, h(n))
\geq \lambda \H^{d-1}(nA)  ] \,>\, 0 \,. $$
\end{thm}

\begin{rem}
We recall the reader that the asymptotic behaviour of $\phi(nA,
h(n))/\H^{d-1}(nA)$ for large $n$ is not known in general. For
straight cylinders, i.e., cylinders of basis $A$ of the form
$\prod_{i=1}^{d-1} [a_i, b_i]\times \{c\}$ with real numbers $a_i$, $b_i$
and $c$, we know
thanks to the works of Kesten \cite{Kesten:flows}, Zhang \cite{Zhang07} and
Rossignol and Th\'eret \cite{RossignolTheret08b} that $\phi(nA, h(n))/\H^{d-1}(nA)$ converges a.s. towards $\nu((0,...,0,1))$
when $n$ goes to infinity, and in this case the upper large deviations of
$\phi(nA,h(n))/\H^{d-1}(nA)$ have been studied by Th\'eret in \cite{TheretUpper}: they are of volume order, and the corresponding large
deviation principle was even proved. For tilted cylinders, we do not know the asymptotic
behaviour of this variable in general, but looking at the trivial case where
$t(e)=1$ for every edge $e$, we can easily see that $\tau(nA,h(n))$ and
$\phi(nA,h(n))$ do not have the same behaviour for large $n$. However, in the
case where $\lim_{n\rightarrow \infty} h(n)/n =0$, we also know that
$\lim_{n\rightarrow \infty} \phi(nA, h(n))/\H^{d-1} (nA) =\nu(\vec{v})$
almost surely under the same hypotheses as for the variable $\tau(nA,h(n))$, so in this case we really study here the upper large
deviations of the variable $\phi(nA, h(n))$.
\end{rem}

\begin{rem}
We were not able to prove a large deviation principle from above for the
variables $\tau$, or $\phi$ in tilted cylinders. The idea used in
\cite{TheretUpper} to prove a large deviation principle for the variable $\phi(nA, h(n))$ in
straight cylinders is the following: we pile up cylinders, and we let a large
amount of flow cross the cylinders one after each other, using the fact
that the top of a cylinder, i.e. the area through which the water goes out
of this cylinder, is exactly the bottom of the cylinder above, i.e. the
area through which the water can go into that cylinder. We cannot use the
same method to prove a large deviation principle for $\tau(nA,h(n))$, even
in straight cylinders, because in this case we cannot glue together the
entire area through which the water goes out of a cylinder with the entire
area through which the water goes into the cylinder above. In the case of
tilted cylinders we even loose the symmetry of the graph
with regard to the hyperplanes spanned by the faces of the cylinder. These
symmetries were of huge importance in the proof of the large deviation
principle from above for $\phi(nA, h(n))/\H^{d-1}(nA)$ in \cite{TheretUpper}. 
\end{rem}

\section[Upper large deviations for $\tau$]{Upper
    large deviations for the rescaled variable $\tau$}

\subsection{Geometric construction}

To study these upper large deviations, we will use the same idea as
in the proof of the strict positivity of the rate function of the large
deviation principle we proved in \cite{TheretUpper} for the variable
$\phi(nA,h(n))$ in straight cylinders. Thus the main tool is the Cram\'er Theorem in $\RR$. We will consider two different scales on
the graph, i.e., cylinders of two different sizes indexed by $n$ and $N$,
with $N$ very large compared to $n$. We
want to divide the cylinder $\cyl (NA, h(N))$ into images of $\cyl(nA, h(n))$
by integer translations, i.e., translations whose vectors have integer coordinates, and to compare the maximal flows through these
cylinders. In fact, we will
first fill $\cyl(NA,h(N))$ with translates of $\cyl(nA,h(n))$ and then move
slightly these translates to obtain integer translates. The problem is that we want
to obtain disjoint small cylinders so that the associated flows are
independent, therefore we need some extra space between the different images of
$\cyl(nA, h(n))$ in order to move them separately and
to obtain disjoint cylinders. Then we add some edges to glue
together the different cutsets in the small cylinders to obtain a
cutset in the big one.

The last remark we have to do before the
beginning of the complete proof is that we may not divide the entire
cylinder $\cyl(NA, h(N))$ into slabs, but a possibly smaller one, $\cyl(NA,
Mh(n))$ with a not too large $M$. Indeed, we will see that the upper
large deviations of $\tau(NA, h(N))$ are related to the behaviour of the
edges of the cylinder that are "not too far" from $NA$, because the cutset
is pinned at the boundary of $NA$ so it cannot explore regions too far away
from $NA$ in $\cyl(NA, h(N))$.

\begin{figure}[ht!]
\centering

\begin{picture}(0,0)%
\includegraphics{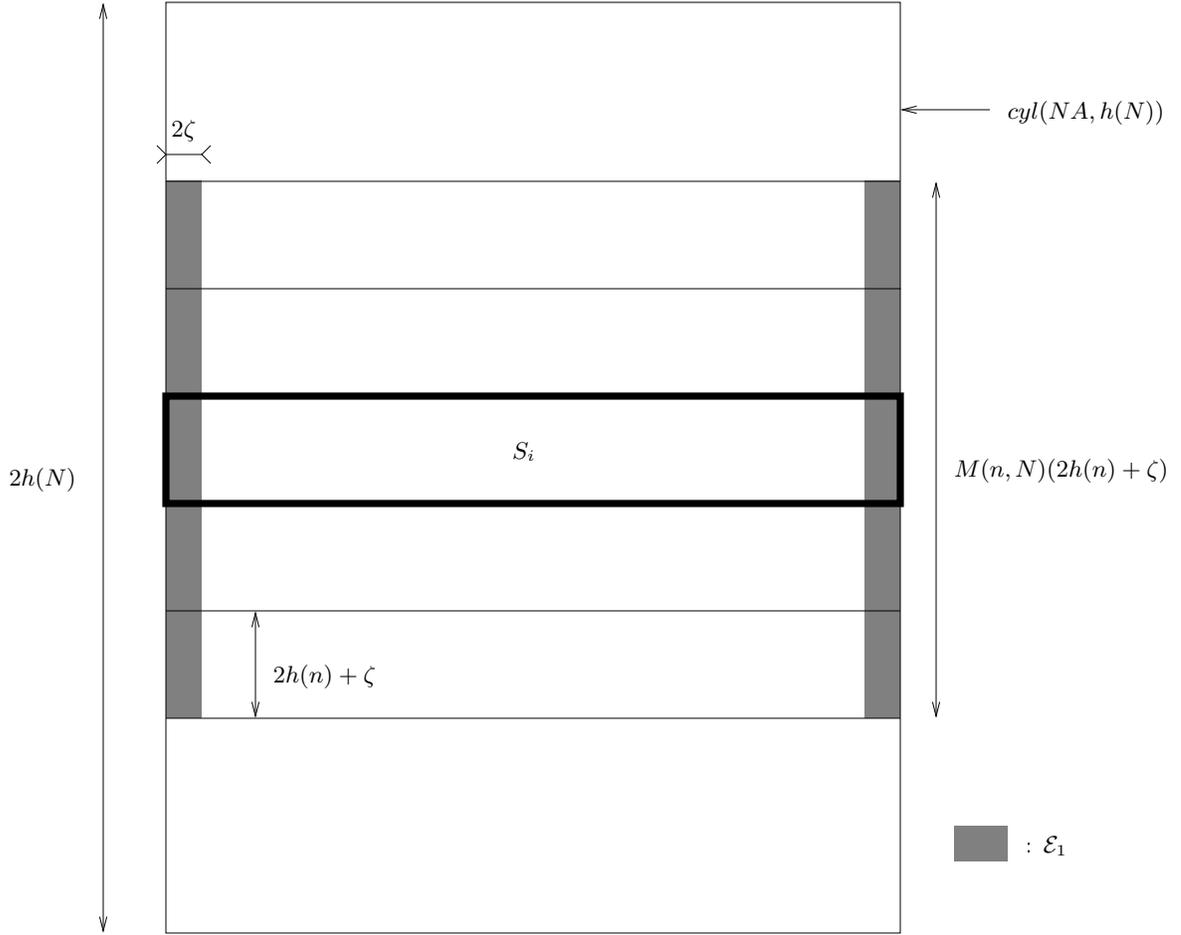}%
\end{picture}%
\setlength{\unitlength}{2960sp}%
\begingroup\makeatletter\ifx\SetFigFont\undefined%
\gdef\SetFigFont#1#2#3#4#5{%
  \reset@font\fontsize{#1}{#2pt}%
  \fontfamily{#3}\fontseries{#4}\fontshape{#5}%
  \selectfont}%
\fi\endgroup%
\begin{picture}(7980,7824)(2536,-6673)
\put(4201,-4561){\makebox(0,0)[lb]{\smash{{\SetFigFont{9}{10.8}{\rmdefault}{\mddefault}{\updefault}{\color[rgb]{0,0,0}$2h(n) + \zeta$}%
}}}}
\put(6301,-2686){\makebox(0,0)[b]{\smash{{\SetFigFont{9}{10.8}{\rmdefault}{\mddefault}{\updefault}{\color[rgb]{0,0,0}$S_i$}%
}}}}
\put(10351,164){\makebox(0,0)[lb]{\smash{{\SetFigFont{9}{10.8}{\rmdefault}{\mddefault}{\updefault}{\color[rgb]{0,0,0}$cyl(NA,h(N))$}%
}}}}
\put(3451, 14){\makebox(0,0)[b]{\smash{{\SetFigFont{9}{10.8}{\rmdefault}{\mddefault}{\updefault}{\color[rgb]{0,0,0}$2\zeta$}%
}}}}
\put(2551,-2911){\makebox(0,0)[rb]{\smash{{\SetFigFont{9}{10.8}{\rmdefault}{\mddefault}{\updefault}{\color[rgb]{0,0,0}$2h(N)$}%
}}}}
\put(10501,-5986){\makebox(0,0)[lb]{\smash{{\SetFigFont{9}{10.8}{\rmdefault}{\mddefault}{\updefault}{\color[rgb]{0,0,0}: $\mathcal{E}_{1}$}%
}}}}
\put(9901,-2836){\makebox(0,0)[lb]{\smash{{\SetFigFont{9}{10.8}{\rmdefault}{\mddefault}{\updefault}{\color[rgb]{0,0,0}$M(n,N) (2h(n) + \zeta)$}%
}}}}
\end{picture}%

\caption{$\cyl(NA,h(N))$ and $S_i$.}
\label{chapitre3devsup1}
\end{figure}

\begin{figure}[ht!]
\centering

\begin{picture}(0,0)%
\includegraphics{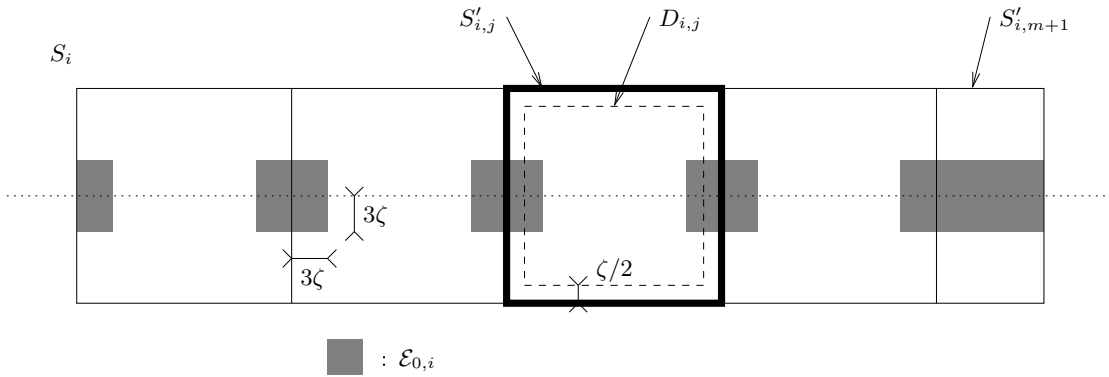}%
\end{picture}%
\setlength{\unitlength}{2960sp}%
\begingroup\makeatletter\ifx\SetFigFont\undefined%
\gdef\SetFigFont#1#2#3#4#5{%
  \reset@font\fontsize{#1}{#2pt}%
  \fontfamily{#3}\fontseries{#4}\fontshape{#5}%
  \selectfont}%
\fi\endgroup%
\begin{picture}(9336,3082)(889,-3662)
\put(5026,-736){\makebox(0,0)[rb]{\smash{{\SetFigFont{9}{10.8}{\rmdefault}{\mddefault}{\updefault}{\color[rgb]{0,0,0}$S'_{i,j}$}%
}}}}
\put(3376,-2911){\makebox(0,0)[lb]{\smash{{\SetFigFont{9}{10.8}{\rmdefault}{\mddefault}{\updefault}{\color[rgb]{0,0,0}$3\zeta$}%
}}}}
\put(3901,-2386){\makebox(0,0)[lb]{\smash{{\SetFigFont{9}{10.8}{\rmdefault}{\mddefault}{\updefault}{\color[rgb]{0,0,0}$3\zeta$}%
}}}}
\put(5851,-2836){\makebox(0,0)[lb]{\smash{{\SetFigFont{9}{10.8}{\rmdefault}{\mddefault}{\updefault}{\color[rgb]{0,0,0}$\zeta /2$}%
}}}}
\put(9226,-736){\makebox(0,0)[lb]{\smash{{\SetFigFont{9}{10.8}{\rmdefault}{\mddefault}{\updefault}{\color[rgb]{0,0,0}$S'_{i,m+1}$}%
}}}}
\put(4051,-3586){\makebox(0,0)[lb]{\smash{{\SetFigFont{9}{10.8}{\rmdefault}{\mddefault}{\updefault}{\color[rgb]{0,0,0}: $\mathcal{E}_{0,i}$}%
}}}}
\put(1276,-1036){\makebox(0,0)[lb]{\smash{{\SetFigFont{9}{10.8}{\rmdefault}{\mddefault}{\updefault}{\color[rgb]{0,0,0}$S_i$}%
}}}}
\put(6376,-736){\makebox(0,0)[lb]{\smash{{\SetFigFont{9}{10.8}{\rmdefault}{\mddefault}{\updefault}{\color[rgb]{0,0,0}$D_{i,j}$}%
}}}}
\end{picture}%

\caption{The slab $S_i$.}
\label{chapitre3devsup2}
\end{figure}

Let $\lambda > \nu(\vec{v})$ and $\varepsilon >0$ such that $\lambda >
\nu(\vec{v}) + 3\varepsilon$. We take an $h$ as in theorem
\ref{chapitre3thmtau}, a large $N$ (we will precise how large it is), and a
smaller $n$. We define $\cyl' (nA,h(n))$ as
$$ \cyl'(nA,h(n)) \,=\, \{x + t\vec{v} \,|\, x\in \hyp(A) \,,\,\,
d(x,nA) \leq \zeta/2 \,\, and \,\, t\in [-h(n)- \zeta/2,h(n)+\zeta/2]  \}\,. $$
We fix an $M=M(n,N)$ such that $M (2 h(n)+\zeta) \leq 2 h(N)$.
We divide $\cyl(NA, M(h(n)+\zeta/2))$ into slabs $S_i$,  $i=1,...,M(n,N)$, of the form
$$ S_i \,=\, \{  x + t\vec{v} \,|\, x\in NA \,,\,\, t\in \mathcal{T}_i \} $$
where
$$ \mathcal{T}_i \,=\, [-M (h(n)+\zeta/2) + (i-1) (2 h(n) +
  \zeta), -M (h(n)+\zeta/2) + i (2 h(n) + \zeta)] $$
(see Figure \ref{chapitre3devsup1}).
By a euclidean division of the dimensions of $S_i$, we divide then each
$S_i$ into $m$ translates of $\cyl'(nA,h(n))$, which we denote by
$S'_{i,j}$, $j=1,...,m$, plus a remaining part $S'_{i,m+1}$.
Here $m$ is smaller than
$\mathcal{M}(n,N) =\lfloor \mathcal{H}^{d-1}(NA) / \mathcal{H}^{d-1}(nA) \rfloor$. Each $S'_{i,j}$ is a
translate of $\cyl'(nA,h(n))$, which contains $\cyl(nA,h(n))$, and so we
denote by $D_{i,j}$ the corresponding translate of $\cyl(nA,h(n))$ by the
same translation ($D_{i,j} \subset S'_{i,j}$). See Figure \ref{chapitre3devsup2} which
illustrates these definitions.

For all $(i,j)$ there exists a vector
$\vec{u}_{i,j}$ in $\mathbb{R}^d$ such that $\|\vec{u}_{i,j}\|_{\infty}
< 1$ and $B_{i,j} = D_{i,j} + \vec{u}_{i,j}$ is the image of
$\cyl(nA,h(n))$ by an integer translation, i.e., a translation whose
vector has integer coordinates; moreover we have $B_{i,j} \subset
S'_{i,j}$, so the $B_{i,j}$ are disjoint. We define $\tau_i =
\tau(S_i,\vec{v})$ and $\tau_{i,j} = \tau(B_{i,j}, \vec{v})$.
We denote by $E_1$ the set of the edges
which belong to $\mathcal{E}_1 \subset \mathbb{R}^d$ defined by
$$ \mathcal{E}_1 \,=\, \{  x + t\vec{v} \,|\, x\in NA \,,\,\, d(x, \partial
(NA)) \leq 2\zeta \,\, and \,\, t\in[-M(h(n)+\zeta/2),M(h(n)+\zeta/2)] \} \,.$$
We denote also by $E_{0,i}$ the set of the edges which belong to
$\mathcal{E}_{0,i} \subset \mathbb{R}^d$ defined by
$$
\mathcal{E}_{0,i} \,=\, \{x+ t\vec{v} \,|\, x \in NA \,,\,\, t\in \mathcal{T}'_i  \}\cap \left( 
\bigcup_{j=1}^{m} \mathcal{V}(\partial S'_{i,j} , 3\zeta) \cup S'_{i,m+1}
\right) \,,
$$
where
$$ \mathcal{T}'_i \,=\, [-h(N) + (i-1/2) (2 h(n) +
  \zeta) - 3 \zeta, -h(N) + (i-1/2) (2 h(n) + \zeta) + 3\zeta] \,.$$
For all $i\in \{1,...,M(n,N)\}$, if we denote by $\F_{i,j}$ a set of edges
that cuts the lower half part from the upper half part of the cylinder
$B_{i,j}$, then $\cup_{j=1}^{m} \F_{i,j} \cup E_{0,i} \cup E_1$ separates
the lower half part from the upper half part of $\cyl(NA, h(N))$. Thus we
obtain that
$$ \forall i\in \{1,...,M(n,N)\}, \qquad \tau(NA, h(N)) \,\leq\, \sum_{j=1}^m
\tau_{i,j} + V(E_1 \cup E_{0,i}) \,, $$
so
\begin{align*} 
\mathbb{P} \Big[ \tau(NA, h(N)) \geq  & \lambda \mathcal{H}^{d-1}(NA)
  \Big]\\ & \,\leq \, \mathbb{P} \left[ \forall i\in \{1,...,M(n,N)\} \,,\,\,
  \sum_{j=1}^m \tau_{i,j} + V(E_1 \cup E_{0,i}) \geq \lambda
  \mathcal{H}^{d-1}(NA)  \right] \\
&\,\leq\, \mathbb{P} \left[ \forall i\in \{1,...,M(n,N)\} \,,\,\, \sum_{j=1}^m \tau_{i,j} \geq
  (\lambda - \varepsilon ) \mathcal{H}^{d-1}(NA) \right]\\
& \qquad  + \mathbb{P}
  \left[ V(E_1 ) \geq
  \varepsilon \mathcal{H}^{d-1}(NA)/2  \right]\\
& \qquad  + \mathbb{P}
  \left[ \exists i\in \{1,...,M(n,N)\} \,,\, V( E_{0,i}) \geq
  \varepsilon \mathcal{H}^{d-1}(NA)/2  \right]\,.\\
\end{align*}
We study the different probabilities appearing here separately.

\noindent $\bullet$ Let 
$$\alpha(N,n) \,=\, \mathbb{P} \left[ \forall i\in
 \{1,...,M(n,N)\} \,,\,\,  \sum_{j=1}^{m} \tau_{i,j} \geq (\lambda -
 \varepsilon ) \mathcal{H}^{d-1}(NA) \right] \,.$$
By independence of the families $(\tau_{i,j}, j=1,...,m)$ for different
 $i$ we have
\begin{align*}
\alpha(N,n) &\,=\, \mathbb{P} \left[ \sum_{j=1}^{m} \tau_{1,j} \geq (\lambda -
  \varepsilon) \mathcal{H}^{d-1}(NA) \right]^{M(n,N)} \\
&\,\leq\, \mathbb{P} \left[ \sum_{j=1}^{\mathcal{M}(n,N,A)} \tau_{1,j} \geq (\lambda -
  \varepsilon) \mathcal{H}^{d-1}(NA) \right]^{M(n,N)} \\
& \,\leq\, \mathbb{P} \left[ \frac{1}{\mathcal{M}(n,N,A)} \sum_{j=1}^{\mathcal{M}(n,N,A)} \frac{\tau_n^{(j)}}{\mathcal{H}^{d-1}(nA)}
  \geq \lambda -\varepsilon \right]^{M(n,N)} \,,\\
\end{align*}
where we remember that 
$$\mathcal{M}(n,N,A) \,=\, \lfloor \mathcal{H}^{d-1}(NA) /
  \mathcal{H}^{d-1}(nA) \rfloor\,, $$
and $(\tau_n^{(j)}, j\in \NN)$ is a family of independent and identically distributed
variables with $\tau_n^{(j)} = \tau(nA, h(n))$ in law. We know that $\mathbb{E}
(\tau(nA, h(n))) / \mathcal{H}^{d-1}(nA)$ converges to $\nu(\vec{v})$ when
$n$ goes to infinity as soon as $\EE[t(e)]<\infty$, so
there exists $n_0$ large enough to have for all $n\geq n_0$ 
$$\frac{\mathbb{E} (\tau(nA, h(n))) }{
\mathcal{H}^{d-1}(nA)} \,\leq\, \nu(\vec{v}) + \varepsilon \,<\, \lambda -
\varepsilon \,.$$
In the three cases presented in Theorem \ref{chapitre3thmtau}, the law of
the capacity of the edges admits at least one exponential moment, and by
an easy comparison between $\tau(nA,h(n))$ and the capacity of a fixed
flat cutset in $\cyl(nA,h(n))$, we obtain that $\tau(nA,h(n))$ admits an
exponential moment. We can then apply the Cram\'er theorem to obtain that for
fixed $n\geq n_0$ and $\lambda$ there exists a constant $c$ (depending on
the law of $\tau(nA,h(n))$, $\lambda$ and $\varepsilon$) such that
$$
\limsup_{N\rightarrow \infty}  \frac{1}{\mathcal{M}(n,N,A)} \log \mathbb{P} \left[ \frac{1}{\mathcal{M}(n,N,A)} \sum_{j=1}^{\mathcal{M}(n,N,A)}
  \frac{\tau_n^{(j)}}{\mathcal{H}^{d-1}(nA)} \geq \lambda - \varepsilon
  \right]  \,\leq\, c \,<\, 0  \,, $$
and so for all $n\geq n_0$ and $\lambda$ there exists a constant $c'$
(depending on the law of $\tau(nA,h(n))$, $\lambda$ and $\varepsilon$) such that
\begin{equation}
\label{chapitre3eqcramer}
 \limsup_{N\rightarrow \infty} \frac{1}{M(n,N)\mathcal{H}^{d-1}(NA)} \log
\alpha(N,n) \,<\, c' <0\,.
\end{equation}

\noindent $\bullet$ To study the two other terms, we can study more
generally the behaviour of
$$ \gamma(n,N) \,=\, \PP\left[ \sum_{i=1}^{l(n,N)} t_i \geq \eps
  \H^{d-1}(NA)/2 \right] \,,$$
where $(t_i, i\in \NN)$ is a family of i.i.d. random variables of common
distribution function $F$. We know that there exists a positive constant $C$ depending on $d$, $A$ and
$\zeta$ such that
\begin{equation}
\label{chapitre3E0}
 \card (E_{0,i}) \,\leq\, C \left( \frac{N^{d-1}}{n} + N^{d-2} n \right)  
\end{equation}
and
\begin{equation}
\label{chapitre3E1}
 \card (E_1) \,\leq\, C N^{d-2} M(n,N) h(n) \,.
\end{equation}
Thus the values of $l(n,N)$ we have to consider are
$$l_0(n,N)\,=\,C (N^{d-1}n^{-1} + N^{d-2}n) \quad \textrm{and} \quad
l_1(n,N)\,=\,CN^{d-2}M(n,N) h(n)\,.$$
The behaviour of the quantity $\gamma(n,N)$ depends on the law of the
capacity of the edges.

\subsection{Bounded capacities}

We suppose that the capacity of the edges is bounded by a constant
$K$. Then as soon as
\begin{equation}
\label{chapitre3condl}
 2 K l(n,N) \,<\, \eps \H^{d-1}(NA)  \,, 
\end{equation}
we know that $\gamma(n,N)=0$. It is obvious that there exists a $n_0$
such that for all fixed $n\geq n_0$, for all large $N$ (how large depending
on $n$), equation (\ref{chapitre3condl}) is satisfied by
$l_0(n,N)$. Moreover, there exists a constant $\kappa (n,A, d ,\zeta,F)$
such that if $M(n,N)\leq \kappa N$, then equation (\ref{chapitre3condl}) is
also satisfied by $l_1(n,N)$ for all $n$. We choose $M(n,N)$ to be as large
as possible according to the condition we have just mentioned, and the
fact that $M(n,N)  \leq h(N)(h(n) +\zeta/2)^{-1}$; we define $\kappa'(n) =
(h(n) +\zeta/2)^{-1}$ and we choose
$$ M(n,N) \,=\, \min(\kappa(n) N, \kappa'(n) h(N) )\,.  $$
Thus, for a fixed $n\geq n_0$,
for all $N$ large enough, we obtain that 
$$ \mathbb{P}
  \left[ V(E_1 ) \geq
  \varepsilon \mathcal{H}^{d-1}(NA)/2  \right] + \mathbb{P}
  \left[ \exists i\in \{1,...,M(n,N)\} \,,\, V( E_{0,i}) \geq
  \varepsilon \mathcal{H}^{d-1}(NA)/2  \right]\,=\,0 $$
and then thanks to equation (\ref{chapitre3eqcramer}) we obtain that
$$ \limsup_{N\rightarrow \infty} \frac{1}{\mathcal{H}^{d-1}(NA) \min (N,
  h(N) )} \log \PP\left[ \frac{\tau(NA, h(N))}{\H^{d-1}(NA)} \geq \lambda
\right] \,<\, 0\,,$$
so equation (\ref{chapitre3casborne}) is proved.

\begin{rem}
The term $\mathcal{H}^{d-1}(nA)\min (n, h(n))$ can seem strange in
(\ref{chapitre3casborne}). It is in fact the right order of the upper large
deviations in the case of bounded capacities. We try here to explain where it comes from. From the point of view of a
minimal cutset, the heuristic is that a cutset in $\cyl(nA, h(n))$
separating the two half cylinders is pinned along the
boundary of $nA$, so it cannot explore domains of $\cyl(nA, h(n))$ that are
too far away from $nA$, i.e., at distance of order larger than $n$. Thus it
is located in a box of volume of order $n^{d-1} \min (n,h(n))$. We
think it is this point of view that gives the best intuitive idea of how
things work, but actually it is very difficult to study the position of a
minimal cutset in the cylinder. From the
point of view of the maximal flow, we can also understand why this term appears. In fact, we can find of the order of $n^{d-1}$
disjoint
paths (i.e., with no common edge) that cross $\cyl(nA, h(n))$ from
its upper half part to its lower half part using only the edges located
at distance smaller
than $K n$ of $nA$ for some constant $K$ (thus all the edges of the box if
$h(n)/n$ is bounded). If $h(n)/n$ is bounded, we can consider paths that cross the
cylinder from its top to its bottom, and if $h(n)\geq n$, we can consider
paths that form a part of a loop around a point of $\partial(nA) $ -
so they join two points of $\cyl(\partial (nA), Kn)$ that are on the same
side of $\cyl(nA, h(n))$ and that are symmetric one to each other by the
reflexion of axis the intersection of $\partial (nA)$ with this side (see
figure \ref{chapitre3minnh} that shows these paths in dimension $2$).
\begin{figure}[ht!]
\centering

\begin{picture}(0,0)%
\includegraphics{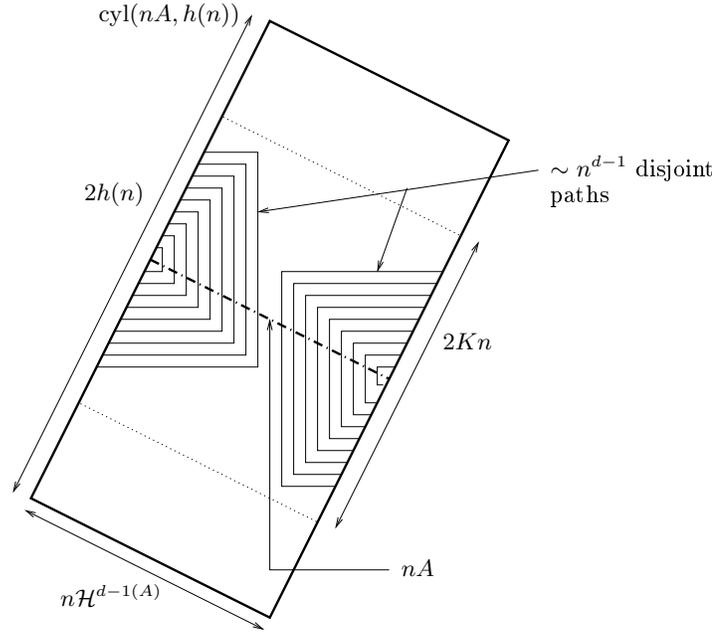}%
\end{picture}%
\setlength{\unitlength}{1973sp}%
\begingroup\makeatletter\ifx\SetFigFont\undefined%
\gdef\SetFigFont#1#2#3#4#5{%
  \reset@font\fontsize{#1}{#2pt}%
  \fontfamily{#3}\fontseries{#4}\fontshape{#5}%
  \selectfont}%
\fi\endgroup%
\begin{picture}(7407,7893)(1264,-7723)
\put(6676,-4186){\makebox(0,0)[lb]{\smash{{\SetFigFont{9}{10.8}{\rmdefault}{\mddefault}{\updefault}{\color[rgb]{0,0,0}$2Kn$}%
}}}}
\put(6151,-7036){\makebox(0,0)[lb]{\smash{{\SetFigFont{9}{10.8}{\rmdefault}{\mddefault}{\updefault}{\color[rgb]{0,0,0}$nA$}%
}}}}
\put(3151,-7411){\makebox(0,0)[rb]{\smash{{\SetFigFont{9}{10.8}{\rmdefault}{\mddefault}{\updefault}{\color[rgb]{0,0,0}$n\H^{d-1(A)}$}%
}}}}
\put(4126,-61){\makebox(0,0)[rb]{\smash{{\SetFigFont{9}{10.8}{\rmdefault}{\mddefault}{\updefault}{\color[rgb]{0,0,0}$\cyl(nA, h(n))$}%
}}}}
\put(2926,-2311){\makebox(0,0)[rb]{\smash{{\SetFigFont{9}{10.8}{\rmdefault}{\mddefault}{\updefault}{\color[rgb]{0,0,0}$2h(n)$}%
}}}}
\put(8026,-2011){\makebox(0,0)[lb]{\smash{{\SetFigFont{9}{10.8}{\rmdefault}{\mddefault}{\updefault}{\color[rgb]{0,0,0}$\sim n^{d-1}$ disjoint}%
}}}}
\put(8026,-2341){\makebox(0,0)[lb]{\smash{{\SetFigFont{9}{10.8}{\rmdefault}{\mddefault}{\updefault}{\color[rgb]{0,0,0}paths}%
}}}}
\end{picture}%

\caption{Disjoint paths near $nA$ in dimension two.}
\label{chapitre3minnh}
\end{figure}
Thus,
if all the edges at distance smaller that $K n$ of $n A$ in the cylinder have a
big capacity, then the variable $\tau(nA, h(n)) /\H^{d-1}(nA)$ is abnormally
big. The number of such edges is of order $n^{d-1} \min(n, h(n))$. We emphasize here the fact that $\phi(nA,
h(n))$ does not have these properties, this is the reason why we expect for
this variable upper large deviations of volume order for all functions $h$.
\end{rem}

\subsection{Capacities of exponential law}

The goal of this short study is to emphasize the fact that the condition of
having one exponential moment for the law of the capacity of the edges is
not sufficient to obtain the speed of decay that we have with bounded
capacities. We will consider a particular law, namely the exponential
law of parameter $1$, and show that we do not have upper deviations of
volume order in this case.

\begin{figure}[ht!]
\centering

\begin{picture}(0,0)%
\includegraphics{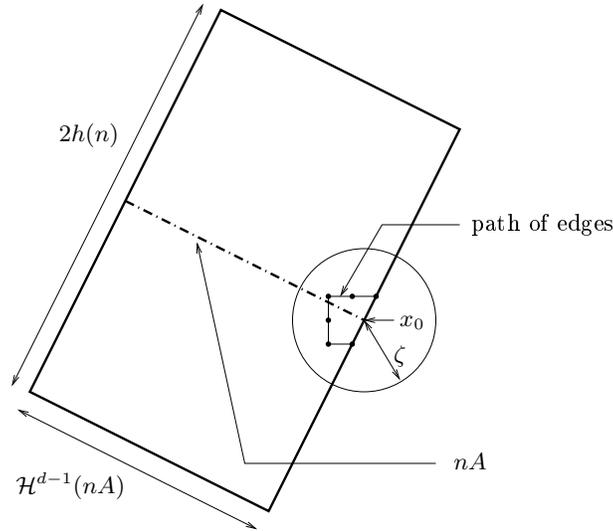}%
\end{picture}%
\setlength{\unitlength}{1973sp}%
\begingroup\makeatletter\ifx\SetFigFont\undefined%
\gdef\SetFigFont#1#2#3#4#5{%
  \reset@font\fontsize{#1}{#2pt}%
  \fontfamily{#3}\fontseries{#4}\fontshape{#5}%
  \selectfont}%
\fi\endgroup%
\begin{picture}(7340,6624)(1564,-7798)
\put(7126,-7036){\makebox(0,0)[lb]{\smash{{\SetFigFont{9}{10.8}{\rmdefault}{\mddefault}{\updefault}{\color[rgb]{0,0,0}$nA$}%
}}}}
\put(3001,-7336){\makebox(0,0)[rb]{\smash{{\SetFigFont{9}{10.8}{\rmdefault}{\mddefault}{\updefault}{\color[rgb]{0,0,0}$\H^{d-1}(nA)$}%
}}}}
\put(2926,-2911){\makebox(0,0)[rb]{\smash{{\SetFigFont{9}{10.8}{\rmdefault}{\mddefault}{\updefault}{\color[rgb]{0,0,0}$2h(n)$}%
}}}}
\put(6451,-5236){\makebox(0,0)[lb]{\smash{{\SetFigFont{9}{10.8}{\rmdefault}{\mddefault}{\updefault}{\color[rgb]{0,0,0}$x_0$}%
}}}}
\put(6376,-5686){\makebox(0,0)[lb]{\smash{{\SetFigFont{9}{10.8}{\rmdefault}{\mddefault}{\updefault}{\color[rgb]{0,0,0}$\zeta$}%
}}}}
\put(7351,-4036){\makebox(0,0)[lb]{\smash{{\SetFigFont{9}{10.8}{\rmdefault}{\mddefault}{\updefault}{\color[rgb]{0,0,0}path of edges}%
}}}}
\end{picture}%

\caption{Path of edges included in a neighbourhood of $x_0$.}
\label{chapitre3xo}
\end{figure}

We suppose that the law of the capacity of the edges is the exponential law
of parameter ~$1$. We know that $\mathbb{E} (\exp (\gamma t))
<\infty$ for all $\gamma <1$. Let $x_0$ be a fixed point of the
boundary $\partial (nA)$. We know that there exists a path from the lower
half cylinder $(nA)_2^{h(n)}$ to the upper half cylinder $(nA)_1^{h(n)}$ in
$\cyl(nA, h(n))$ that is included in the neighbourhood of $x_0$ of diameter
$\zeta \geq 2d$ for the euclidean distance, as soon as $n \geq n_0(d,A,h)$,
where $n_0(d,A,h)$ is the infimum of the $n$ such that all the sidelengths
of the cylinder $\cyl(nA, h(n))$ are larger than $\zeta$ (see figure
\ref{chapitre3xo}). 
Thus for all $n\geq n_0$, every set of edges that
cuts the upper half cylinder $(nA)_1^{h(n)}$ from the lower half cylinder
$(nA)_2^{h(n)}$ in $\cyl (nA, h(n))$ must contain
one of the edges of this neighbourhood of $x_0$. The number of such
edges is at most $K(d,\zeta)$, where $K$ is a constant depending only on
$d$ and $\zeta$.
Thus the probability that all of them have a capacity
bigger than $\lambda \mathcal{H}^{d-1}(nA)$ for a $\lambda > \nu(\vec{v})$
is greater than $\exp (-K \lambda \mathcal{H}^{d-1}(nA) )$. We obtain that
for all $n\geq n_0(d,A,h)$,
$$ \PP \left[ \tau (nA, h(n))\geq \lambda \mathcal{H}^{d-1}(nA) \right]
  \,\geq\,\exp (-K \lambda \mathcal{H}^{d-1}(nA) ) \,.$$

\subsection[Capacities with exp. moment of all orders]{Capacities with exponential moments of all orders}

We suppose that the capacity of the edges admits exponential
moments of all order, i.e., for all $\theta >0$ we have $\EE(\exp (\theta t(e)))
<\infty$. Then by a simple application of the Chebyshev inequality, we
obtain that
\begin{equation}
\label{chapitre3eqtcheb}
 \gamma(n,N) \,\leq\, \exp \left[ -\H^{d-1}(NA) \left(  \frac{\theta \eps
      }{2 } - \frac{l(n,N)\log \EE (\exp(\theta t(e)))}{\H^{d-1}(NA)} \right)
\right] \,.
\end{equation}
We want to be able to choose the term 
$$ \frac{\theta \eps
      }{2 } - \frac{l(n,N)\log \EE (\exp(\theta t(e)))}{\H^{d-1}(NA)}$$
as big as we want. 
For a fixed $R>0$, we can take $\theta >0$ large enough to have $\theta
\varepsilon \geq 4R$. If there exists $n_1$ such that for a fixed $n\geq n_1$,
for all $N$ sufficiently large (how large depends on $n$), we have
\begin{equation}
\label{chapitre3condl2}
\frac{l(n,N)}{\H^{d-1}(NA)} \log \EE(e^{\theta t(e)}) \,\leq\, R\,,
\end{equation}
then for a fixed $n\geq n_1$, for all large $N$, we would obtain
$$ \gamma(n,N) \,\leq\,  \exp \left(-R \mathcal{H}^{d-1}(NA)\right) \,. $$
We consider now the values of $l_0(n,N)$ and $l_1(n,N)$.
Looking at $l_1(n,N)$, we realize that we have to impose a
condition on $M(n,N)$. Considering the result we want to prove, we can
choose $M(n,N)$ satisfying, for each fixed $n$,
$$ \lim_{N\rightarrow \infty} \frac{M(n,N)}{N} \,=\, 0  \qquad \textrm{and} \qquad
\lim_{N\rightarrow \infty} M(n,N) \,=\, +\infty \,. $$
If $h(N)/N$ does not converge towards $0$, we thus consider a small
cylinder inside the cylinder $\cyl(NA, h(N))$, but we impose that its
height goes to
infinity with $N$. Under this hypothesis, we obtain that for all $R$, for
every fixed $n$, for all large $N$, condition (\ref{chapitre3condl2}) is
satisfied by $l_1(n,N)$. Thus, for all fixed $n$, thanks to
(\ref{chapitre3E1}), we obtain that
\begin{equation}
\label{chapitre3mtsb}
\limsup_{N\rightarrow \infty} \frac{1}{\H^{d-1}(NA)} \log \PP  \left[  V(
  E_{1}) \geq   \varepsilon \mathcal{H}^{d-1}(NA)/2  \right] \,=\, -\infty\,.
\end{equation}
For all $R$, we can find a $n_1$ such that for all $n\geq n_1$, for all
large $N$, the condition (\ref{chapitre3condl2}) is satisfied by
$l_0(n,N)$. Since our choice of $M(n,N)$ implies that 
$$ \lim_{N\rightarrow \infty} \frac{\log M(n,N)}{\H^{d-1}(NA)} \,=\,0\,,    $$ 
thanks to
(\ref{chapitre3E0}), we obtain that for all fixed $n\geq n_1$,
\begin{equation}
\label{chapitre3mtsa}
\limsup_{N\rightarrow \infty} \frac{1}{\H^{d-1}(NA)} \log \PP  \left[ \exists i\in \{1,...,M(n,N)\} \,,\, V( E_{0,i}) \geq
  \varepsilon \mathcal{H}^{d-1}(NA)/2  \right] \,=\, -\infty\,.
\end{equation}
Combining (\ref{chapitre3mtsb}), (\ref{chapitre3mtsa}) and
(\ref{chapitre3eqcramer}), since $\lim_{N\rightarrow \infty} M(n,N)
=+\infty$, we have proved (\ref{chapitre3castsmts}). This ends the proof of
Theorem \ref{chapitre3thmtau}.

\begin{rem}
This result is used in \cite{RossignolTheret08b} in the proof of the lower large
deviation principle for the variable $\tau(nA,h(n))$.
\end{rem}

%%%%%%%%%%%%%%%%%%%%%%%%%%%%%%%%%%%%%%%%%%%%%%%%%%%%%%%%%%%%%%%%%%%%%%%%%%%%%%%

\section[Upper large deviations for $\phi$]{Partial result concerning the upper large
    deviations for $\phi$ through a tilted cylinder}

We have already written the main part of the proof of Theorem
\ref{chapitre3thmphi} in the previous section. We keep all the notations
introduced previously. The proof of Theorem \ref{chapitre3thmtau} was based
on the following inequality:
$$ \forall i\in \{1,...,M(n,N)\}, \qquad \tau(NA, h(N)) \,\leq\, \sum_{j=1}^m
\tau_{i,j} + V(E_1 \cup E_{0,i}) \,. $$
We recall that this inequality was obtained by noticing that if $\F_{i,j}$
is a cutset that separates the upper half part from the lower half part of
$B_{i,j}$, then $\cup_{j=1}^m \F_{i,j} \cup E_{0,i} \cup E_1$ separates the
upper half part from the lower half part of $\cyl(NA, h(N))$. Here we want
to construct a cutset that separates the bottom from the top of $\cyl(NA,
h(N))$. We have no need to add the set of edges $E_1$ in this context
because we do not need to obtain a cutset that is pinned at $\partial
(NA)$. Thus for all $i$, $\cup_{j=1}^m \F_{i,j} \cup E_{0,i}$ cuts the top
from the bottom of $\cyl(NA, h(N))$, and then we have
$$\forall i\in \{1,...,M(n,N)\} \,,\qquad \phi(NA, h(N)) \,\leq\,
\sum_{j=1}^m \tau_{i,j} +V(E_{0,i})\,.  $$
We obtain that for a fixed $\lambda > \nu(\vec{v})$, and $\eps$ such
that $\lambda \geq \nu(\vec{v}) + 3\eps$, we have by independence
\begin{align*}
\PP [ & \phi(NA, h(N) )  \geq \lambda \H^{d-1}(NA) ]\\
& \,\leq\,\PP \left[\bigcap_{i=1}^{M(n,N)} \left\{ \sum_{j=1}^m \tau_{i,j} +V(E_{0,i})
    \geq \lambda \H^{d-1}(NA) \right\} \right]\\
&  \,\leq\, \prod_{i=1}^{M(n,N)}\left( \PP \left[ \sum_{j=1}^m \tau_{i,j}
    \geq (\lambda - \eps) \H^{d-1}(NA)  \right] + \PP \left[V(E_i)\geq \eps \H^{d-1}(NA)
  \right] \right)\,. 
\end{align*}
We consider here the maximal $M(n,N)$, i.e., 
$$M(n,N) \,=\, \left\lfloor \frac{h(N)}{h(n) + \zeta/2} \right\rfloor \,.$$
Indeed, we do not need to make any restriction on $M(n,N)$ because we do not
have to consider the set of edges $E_1$ whose cardinality depends on
$M(n,N)$.

From now on we suppose that the capacity of the edges admits an exponential
moment. Thanks to the application of the Cram\'er theorem we have already
done to obtain (\ref{chapitre3eqcramer}), we know that for all $n\geq n_0$
there exists a positive $c'$ (depending on the law of $\tau(nA,h(n))$,
$\lambda$ and $\varepsilon$) such that
\begin{equation}
\label{chapitre3resphia}
\limsup_{N\rightarrow \infty} \frac{1}{\H^{d-1} (NA) } \log \PP \left[
  \sum_{j=1}^{m} \tau_{i,j} \geq (\lambda - \varepsilon)\H^{d-1} (NA)
\right] \,\leq\, c' \,<\,0\,.
\end{equation}

On the other hand, let $\gamma >0$ be such that $\EE(\exp (\gamma t(e)))
<\infty$. Thanks to equation (\ref{chapitre3eqtcheb}), obtained by the
Chebyshev inequality, and (\ref{chapitre3E0}), we have for this fixed $\gamma$:
\begin{align*}
 \PP [V(E_{0,i}) \geq \eps \H^{d-1}(NA)]& \,\leq\, \PP \left[
  \sum_{i=1}^{l_0(n,N)} t_i \geq \eps \H^{d-1}(NA)\right]  \\
& \,\leq\,  \exp \left[ -\H^{d-1}(NA) \left(  \frac{\gamma \eps
      }{2 } - \frac{l_0(n,N)\log \EE (\exp(\gamma t(e)))}{\H^{d-1}(NA)} \right)
\right] \,.
\end{align*}
Since $l_0(n,N) \leq C(N^{d-1} n^{-1} + N^{d-2}n)$, we know that there
exists $n_1$ such that for all $n\geq n_1$, for all $N$ large enough
(how large depending on $n$), we have
$$ \frac{l_0(n,N) \log \EE(\exp (\gamma t(e)))}{\H^{d-1}(NA)} \,\leq\,
\frac{\gamma \eps}{ 4} \,, $$
and then
\begin{equation}
\label{chapitre3resphib}
 \PP [V(E_{0,i}) \geq \eps \H^{d-1}(NA)] \,\leq\, \exp \left( -\H^{d-1}(NA)
 \frac{\gamma \eps}{ 4} \right)\,.
\end{equation}
Combining equations (\ref{chapitre3resphia}) and (\ref{chapitre3resphib}),
since $M(n,N)$ is proportional to $h(N)$ for a fixed $n$, Theorem
\ref{chapitre3thmphi} is proved.\\

{\bf Acknowledgements}

The author want to thank Rapha\"el Cerf and Rapha\"el Rossignol, who asked
the questions that have motivated this work.

%%%%%%%%%%%%%%%%%%%%%%%%%%%%%%%%%%%%%%
\def\cprime{$'$}

%\bibliographystyle{plain}
%\bibliography{biblio}

\end{document}